\documentclass[a4paper,10pt]{amsart}
\usepackage[english]{certus}
 \usepackage{bm}

\newcommand{\tens}{\otimes}

\newcommand{\Tor}{{\operatorname{Tor}}}
\newcommand{\Ext}{{\operatorname{Ext}}}

\newcommand{\twoone}{\operatorname{II}_1}

\theoremstyle{definition}

\newcommand{\Proj}{{\operatorname{Proj}}}

\title{On the first continuous $L^2$-cohomology of free group factors}
\author{Vadim Alekseev}

\address{Vadim Alekseev,
Mathematisches Institut,
Georg-Au\-gust-Uni\-versi\-t{\"a}t G{\"o}t\-ting\-en,
Bunsenstra{\ss}e 3-5,
D-37073 G{\"o}ttingen, 
Germany.}
\email{alekseev@uni-math.gwdg.de}
\urladdr{www.uni-math.gwdg.de/alekseev}


\begin{document}

 \begin{abstract}
We prove that the first continuous $L^2$-cohomology of free group factors vanishes. This answers a question by Andreas Thom regarding continuity properties of free difference quotients and shows that one can not distinguish free group factors by means of first continuous $L^2$-Betti number.
 \end{abstract}
\maketitle
 
\section{Introduction}
Introduced by topologists \cite{atiyah-l2}, $L^2$-Betti numbers have been generalized to various contexts like groups, groupoids etc. Alain Connes and Dimitri Shlyakhtenko \cite{conshl2005l2} introduced $L^2$-Betti numbers for subalgebras of finite von Neumann algebras, with the purpose to obtain a suitable notion for arbitrary II$_1$-factors and in the hope to get a nice homological invariant for them. Unfortunately, as of now there are only very few concrete calculations of them. The most advanced computational result so far is due to Andreas Thom \cite{thom2008l2} who proved that the $L^2$-Betti numbers vanish for von Neumann algebras with diffuse center. To allow more computable examples, he also introduced a continuous version of the first $L^2$-Betti number \cite{thom2008l2} which turns out to be much more manageable than its algebraic counterpart. The first continuous $L^2$-Betti number is defined as the von Neumann dimension of the first continuous Hochschild cohomology of the von Neumann 
algebra $M$ 
with values in the algebra of operators affiliated with $M\htens M^\op$. The word `continuous' here means that we restrict attention to derivations which are continuous from the norm topology on $M$ to the measure topology on the affiliated operators.

So far only vanishing results were obtained about the first continuous $L^2$-Betti number: it has been shown to vanish for II$_1$-factors with Cartan subalgebras, non-prime II$_1$-factors \cite{thom2008l2} as well as for II$_1$-factors with property (T), property $\Gamma$ and finitely generated II$_1$ factors with nontrivial fundamental group \cite{alekseevkyed2012-first-l2}. The last result is due to a compression formula for the first continuous $L^2$-Betti number \cite[Theorem 4.10]{alekseevkyed2012-first-l2}.

The hope placed upon $L^2$-Betti numbers for group von Neumann algebras was to be able to connect them with $L^2$-Betti numbers of groups, thus obtaining a powerful invariant which would be able to distinguish free group factors, thus solving a long-standing problem in operator algebras. In fact, the attempt to do this can be formulated in a very concrete way using generators of the $L^2$-cohomology of the group ring $\mb C\mb F_n$ of the free group or some other subalgebras of $L\mb F_n$ generated by free elements. One possible choice of generators is to consider the so-called Voiculescu's free difference quotients \cite{voiculescu-entropy-v}. Andreas Thom posed a natural question in \cite{thom2008l2}, whether these derivations possess continuous extensions to operators from $L\mb F_n$ to $\ms U(L\mb F_n\htens L\mb F_n^\op)$; a positive answer to this question would solve the free factor isomorphism problem.

In the present paper we answer this question in the negative; in fact, we show that the first continuous $L^2$-cohomology of free group factors vanishes; in particular, they can not be distinguished by this invariant. This also suggests that the invariant might be altogether trivial, i.e. that the first continuous $L^2$-cohomology might in fact vanish for all II$_1$-factors.

The result is established in several steps. First, we focus on the free group with three generators $\mb F_3$ and show that the canonical derivations which ``derive in direction of a free generator'' cannot be extended to the group von Neumann algebra. This is shown by analyzing their values on some specific elements for which the spectrum of the resulting operators can be calculated using free probability theory. To derive the vanishing of the whole continuous cohomology, we have to use certain automorphisms of the free group factors. Hereby we make use of certain weak mixing properties relative to a subalgebra; intuitively speaking, we are using the fact that there are enough automorphisms to move our derivations around; thus, the existence of one continuous non-inner derivation would automatically guarantee that all derivations of $\mb C \mb F_3$ are extendable, which yields a contradiction. Finally, we make use of the compression formula to extend the result from a single free group factor to all of them.

The author thanks Thomas Schick and Andreas Thom for helpful discussions and useful suggestions.

\section{Preparatory results}
In this section we set up the notation and briefly recapitulate the theory of non-commutative integration and the theory of $L^2$-Betti numbers for von Neumann algebras.

\subsection{Notation}
We consider finite von Neumann algebras $M$, $N$ etc. with separable preduals. We always endow them with a fixed faithful normal tracial state (usually denoted by $\tau$) and consider them in the corresponding GNS representation $L^2(N,\tau)$. If $(N,\tau)$ is a finite von Neumann algebra, then there is an induced a faithful normal tracial state on the von Neumann algebraic tensor product $N\overline{\otimes} N^{\textup{op}}$ of $N$ with its opposite algebra; abusing notation slightly, we will still denote it by $\tau$. We let $\ms U(N)$ be the algebra of closed densely defined operators on $L^2(N,\tau)$ affiliated with $N$. We equip $\ms U(N)$ with the \emph{measure topology}, defined by the following two-parameter family of zero neighbourhoods:
\beqn
N(\eps,\delta) = \{a\in \ms U(N)\,|\, \exists p\in \mathrm{Proj}(N): \norm{ap} < \eps,\; \tau(p^\perp) < \delta\}, \quad \eps,\delta>0.
\eeqn
With this topology, $\ms U(N)$ is a complete \cite[Theorem IX.2.5]{takesaki2003theory2}  metrizable \cite[Theorem 1.24]{rudin-funct-an}  topological vector space  and the multiplication map
\beqn
(a,b)\mapsto ab \colon \ms U(N)\times \ms U(N) \to \ms U(N) 
\eeqn
is uniformly continuous when restricted to products of bounded subsets \cite[Theorem 1]{nelson-nc-integration}. Convergence with respect to the measure topology is also referred to as \emph{convergence in measure} and denoted by $\xrightarrow{m}$. If $\xi\in\ms U(N)$ and $p\in N$ is its source projection, we denote $\rk \xi:= \tau(p)$. Of course, we also have $\rk\xi = \tau(q)$, where $q$ is the target projection of $\xi$.

Here and in the sequel $\odot$ denotes the algebric tensor product over $\mb C$. We freely identify $M$-$M$-bimodules with $M\odot M^\op$-modules. For $N=M\htens M^\op$ we equip $\ms U(M\htens M^\op)$ with the $M$-$M$-bimodule structure
\[
m\cdot \xi:= (m\tens1^\op)\xi \ \text{ and } \ \xi\cdot m:= (1\tens m^\op)\xi \ \text{ for } m,m \in M \text{ and } \xi \in \ms U. 
\]
All $M$-$M$-sub-bimodules of $\ms U(M\htens M^\op)$ inherit this bimodule structure.

Let $\Gamma$ be a discrete group with its natural left action $\lambda$ on itself. Its von Neumann algebra $L\Gamma = \lambda(\Gamma)''\subset \mb B(\ell^2(\Gamma))$ is equipped with the natural faithful normal tracial state $\tau(\cdot):=\ip{\cdot \delta_e,\delta_e}$. Notice that the GNS representation of $L\Gamma$ with respect to $\tau$ coincides with $\ell^2(\Gamma)$.

Let $\mc A$ be an algebra and $\mc X$ an $\mc A$-$\mc A$-bimodule. We recall that a linear map $\delta\colon \mc A\to \mc X$ is called a \emph{derivation} if it satisfies
\[
\delta(ab)=a\cdot \delta(b) + \delta(a)\cdot b \ \text{ for all } a,b\in \mc A,
\]
and that a derivation is called \emph{inner} if there exists a vector $\xi \in \mc X$ such that
\[
\delta(a)=a\cdot \xi-\xi\cdot a \text{ for all } a\in \mc A.
\]

Consider a free group $\mb F$ generated by a set $S$. For $s\in S$, a derivation
\[
 \partial_s\colon \mb C\mb F\to L\mb F\htens L\mb F^\op
\]
is defined uniquely by the properties
\[
 \partial_s(s) = s\otimes 1^\op,\quad \partial_s(s') = 0,\quad s\neq s'\in S.
\]


\subsection{Some properties of convergence in measure} 
Here we prove several lemmas which will help us to analyse convergence in measure of some particular operators. The following lemma tells us that we can analyse convergence in measure ``locally''.

\begin{Lemma}\label{lemma:measure-conv-local-crit}
 Let $\{a_n\}_{n=1}^\infty\subset \ms U(N)$ be a sequence. Then $a_n\xrightarrow{m} 0$ if and only if
 \beq\label{eq:conv-in-measure-to-zero}
  \forall \eps > 0\quad \tau\left(\chi_{[\eps^2,+\infty)}(a_n^* a_n)\right) \to 0,\quad n\to \infty.
 \eeq
\end{Lemma}
\begin{proof}
 Suppose that \eqref{eq:conv-in-measure-to-zero} is satisfied. For every $\eps > 0$ set $p_n := 1 - \chi_{[\eps^2,+\infty)}(a^*a)$. We obtain
 \[
 \norm{a_n p_n}^2 = \norm{p_na_n^*a_n p_n} \leqslant \eps^2 
 \]
 and
 \[
  \tau(p_n) \to 1,\quad n\to\infty.
 \]
 Thus, $a_n\xrightarrow{m} 0$.
 
 On the other hand, let $a_n \xrightarrow{m} 0$ and suppose that \eqref{eq:conv-in-measure-to-zero} is false. Then, extracting a subsequence if needed, we may assume that 
 \[
  \tau\left(\chi_{[\eps^2,+\infty)}(a_n^*a_n)\right) \geqslant \delta,\quad n\geqslant N_0.
 \]
Suppose now that there exists a sequence of projections $\{q_n\}_{n=1}^\infty$ such that
\[
 \norm{a_nq_n}\to 0,\quad \tau(q_n) \to 1.
\]
Then for sufficiently big $n$ we get $\tau(q_n) > 1 - \delta/2$. Putting $r_n:=\chi_{[\eps^2,+\infty)}(a_n^*a_n)$, we obtain
\[
\norm{a_nq_n}^2 = \norm{q_na_n^*a_n q_n}\geqslant \norm{r'_n a_n^*a_n r'_n} \geqslant \eps^2
\]
for $r'_n:= r_n\wedge q_n$, which satisfies $\tau(r'_n)\geqslant \delta/2$, obtaining a contradiction.
\end{proof}

\begin{Lemma}\label{lemma:no-convergence-to-zero}
 Let $a_n\in \ms U(N)$ be a sequence of selfadjoint elements such that
 \[
  \forall \eps > 0\quad \tau\left(\chi_{[-\eps,\eps]}(a_n)\right) \to 0,\quad n\to\infty.
 \]
Then for every nonzero projection $p\in N$
\[
a_n p \not \xrightarrow{m} 0.
\]
\end{Lemma}
\begin{proof}
 From the assumption of the lemma it immediately follows that
 \[
  \forall \eps > 0\quad \tau\left(\chi_{[\eps^2,+\infty)}(a_n^*a_n)\right) \to 1.
 \]
Setting $r_n:=\chi_{[\eps^2,+\infty)}(a_n^*a_n)$, we get $\tau(r_n) \to 1,\quad n\to\infty$. Now, if $\tau(p) = \delta > 0$, then there exists an $N_0$ such that $\tau(r_n) > 1-\delta/2,\; n\geqslant N_0$. Thus for such $n$ we obtain that the projections $q_n:= r_n\wedge p$ satisfy $\tau(q_n)\geqslant \delta/2$. It follows that
\[
q_na_n^*a_n q_n \geqslant \eps^2q_n,
\]
and hence by Lemma \ref{lemma:measure-conv-local-crit} $q_na_n^*a_n q_n \not\xrightarrow{m} 0$. On the other hand, if $a_np\xrightarrow{m} 0$, then by boundedness of $q_n$
\[
 q_na_n^*a_n q_n = q_npa_n^*a_np q_n \xrightarrow{m} 0.
\]
This proves the lemma.
\end{proof}

\subsection{\texorpdfstring{$\bm{L}^{\bm2}$-Betti numbers for groups and operator algebras}{L2-Betti numbers for groups and operator algebras}}
In \cite{conshl2005l2} Connes and Shlyakhtenko introduced $L^2$-Betti numbers in the general setting of tracial $*$-algebras; if $M$ is a finite von Neumann algebra and $\mc A\subset M$ is any weakly dense unital $\ast$-subalgebra its $L^2$-Betti numbers are defined as
\[
\beta_p^{(2)}(\mc A, \tau)=\dim_{M\htens M^\op} \Tor_p^{\mc A \odot \mc A^\op }(M\htens M^\op, \mc A).
\]
Here the dimension function $\dim_{M\htens M^\op}(-)$ is L{\"u}ck's extended von Neumann dimension~\cite[Chapter 6]{luck2002book}.
This definition is inspired by the well-known correspondence between representations of groups and bimodules over finite von Neumann algebras, and it extends the classical theory by means of the formula $\beta_p^{(2)}(\Gamma)=\beta_p^{(2)}(\mathbb{C}\Gamma,\tau)$ whenever $\Gamma$ is a discrete countable group. In \cite{thom2008l2} it is shown that the $L^2$-Betti numbers also allow the following cohomological description:
\[
\beta_p^{(2)}(\mc A,\tau)=\dim_{M\htens M^\op}\Ext_{\mc A \odot \mc A}^p(\mc A, \ms U),
\]
where $\ms U = \ms U(M\htens M^\op)$ denotes the algebra of operators affiliated with $M\htens M^\op$. It is a classical fact \cite[1.5.8]{loday1998cyclic} that the Ext-groups above are isomorphic to the Hochschild cohomology groups of $\mc A$ with coefficients in  $\ms U$, where the latter is considered as an $\mc A$-bimodule with respect to the actions
\[
a\cdot \xi:= (a\tens1^\op)\xi \ \text{ and } \ \xi\cdot b:= (1\tens b^\op)\xi \ \text{ for } a,b \in \mc A \text{ and } \xi \in \ms U. 
\]
In particular, the first $L^2$-Betti number can be computed as the dimension of the right $M\htens M^\op$-module
\[
H^1(\mc A, \ms U)=\frac{\Der(\mc A, \ms U)}{\text{Inn}(\mc A,\ms U)}. 
\]
Here $\Der(\mc A,\ms U) $ denotes the space of derivations from $\mc A$ to $\ms U$ and $\text{Inn}(A,\ms U)$ denotes the space of inner derivations. 
These purely algebraically defined $L^2$-Betti numbers have turned out extremely difficult to compute in the case when $\mc A$ is $M$ itself. Besides finite-dimensional algebras, the only computational result known in this direction vanishing for von Neumann algebras with diffuse centre  (see \cite[Corollary 3.5]{conshl2005l2} and \cite[Theorem 2.2]{thom2008l2}). It is therefore natural to consider variations of the definitions above that take into account the topological nature of $M$. The continuous version of the first $L^2$-cohomology module was introduced by Andreas Thom in \cite{thom2008l2}, where one restricts attention to those derivations $\delta\colon \mc A \to \ms U$ which are closable from the norm topology to the measure topology. Note that when $\mc A$ is norm closed these are exactly the derivations that are norm--measure topology continuous by the closed graph theorem. We denote the space of closable derivations by $\Der_c(\mc A, \ms U)$, the  continuous cohomology by $H^1_c(\mc A, \ms U)$ 
and by $\eta_1^{(2)}(\mc A,\tau)$ the corresponding continuous $L^2$-Betti numbers; i.e.~ 
\[
\eta_1^{(2)}(\mc A,\tau)=\dim_{M\htens M^\op} H^1_c(\mc A,\ms U).
\]
Notice that by continuity of multiplication on $\ms U$, $\Der_c(\mc A, \ms U)$ is naturally a right $\ms U$-module.

The first continuous $L^2$-Betti number satisfies the following compression formula analogous to \cite[Theorem 2.4]{conshl2005l2}:
\begin{Thm}[{\cite[Theorem 4.10]{alekseevkyed2012-first-l2}}]\label{thm:compression-formula}
 Let $M$ be a $\twoone$-factor with trace-state $\tau$ and let $p\in M$ be a non-zero projection. Then $\eta_1^{(2)}(pMp, \tau_p)=\frac{1}{\tau(p)^2}\eta_1^{(2)}(M,\tau)$.
\end{Thm}

Although the extended von Neumann dimension is generally not faithful, enlarging the coefficients from $M\htens M^\op$ to $\ms U$  has the effect that $\beta_1^{(2)}(\mc A,\tau)=0$ if and only if $H^1(\mc A, \ms U)$ vanishes \cite[Corollary 3.3 and Theorem 3.5]{thom2008l2}. In particular, in order to prove that $\beta_1^{(2)}(\mc A, \tau)=0$ one has to prove that every derivation from $\mc A$ into $\ms U$ is inner. An analogous statement holds for continuous $L^2$-cohomology:
\begin{Prop}[{\cite[Proposition 4.3]{alekseevkyed2012-first-l2}}]\label{prop:vanishing-cohomology-vs-vanishing-betti}
 Let $M$ be a finite von Neumann algebra. We have $\eta_1^{(2)}(M)=0$ if and only if $H^1_c(M,\ms U)\cong 0$.
\end{Prop}

\subsection{Continuity properties of derivations}
We will be interested in continuity of certain derivations. To understand it, we recollect some useful notions and properties here.

Let $\mc A\subset M$ be a weakly dense $\ast$-subalgebra and $\delta\colon \mc A\to \ms U$ be a derivation. Let
\[
 P_\delta:=\{p\in\Proj(M\htens M^\op)\,|\, \delta\cdot p\text{ is continuous}\}.
\]
\begin{Lemma}\label{lemma:complete-lattice}
 $P_\delta$ is a complete sublattice of $\Proj(M\htens M^\op)$.
\end{Lemma}
\begin{proof}
 If $\xi \in\ms U$, the right $\ms U$-submodule generated by $\delta\cdot \xi$ contains the derivation $\delta\cdot t(\xi)$, where $t(\xi)$ is the target projection of $\xi$. If now $p_1,p_2\in P_\delta$, then $p_1\wedge p_2\in P_\delta$ for obvious reasons, and $p_1\vee p_2\in P_\delta$ because $\delta\cdot(p_1+p_2)$ is continuous and $p_1\vee p_2 = t(p_1+p_2)$. If now $\{p_i\}_{i\in I}$ is an orthogonal family of projections with sum $p$, then
 \[
  \sum_{i\in I} \tau(p_i) = \tau(p)\leqslant 1,
 \]
and therefore the series
\[
 \sum_{i\in I} \delta(x)\cdot p_i
\]
converges uniformly in measure to $\delta(x)\cdot p$. As sums of uniformly convergent series of continuous maps are continuous, $\delta\cdot p$ is continuous.
\end{proof}

\begin{Def}
 We call the unique supremum of $P_\delta$ the \emph{continuity projection} of $\delta$.
\end{Def}

Let $\delta\colon \mc A \to \ms U$ be a derivation and $\sigma\in\Aut (M)$ be an automorphism; by slight abuse of notation, we still denote by $\sigma$ the induced automorphism $\sigma\otimes \sigma^\op\in\Aut (\ms U)$. The map
\[
 \delta^\sigma:=\sigma^{-1}\circ \delta\circ \sigma
\]
is then a derivation $M \to \ms U$. If $\delta$ is continuous, then $\delta^\sigma$ is obviously continuous as well. Notice that for $\xi\in \ms U$ we have
\[
 (\delta\cdot \xi)^\sigma = \delta^\sigma \cdot \sigma^{-1}(\xi).
\]

The following observation is easy, but very useful.

\begin{Lemma}\label{lemma:invariant-projection}
If $\Sigma\subset \Aut(M)$ is a subgroup of automorphisms of $M$ and $\delta^\sigma = \delta$ for $\sigma\in \Sigma$, then the continuity projection of $\delta$ is $\Sigma$-invariant.
\end{Lemma}
\begin{proof}
Let $p$ be the continuity projection of $\delta$. Then $\delta\cdot p$ is continuous, and therefore
\[
  (\delta\cdot p)^\sigma = \delta^\sigma \cdot \sigma^{-1}(p).
\]
is continuous as well. If for some $\sigma\in \Sigma$ we have $\sigma^{-1}(p) \neq p$, then $\sigma^{-1}(p) \vee p \in P_\delta$ has $p$ as a proper subprojection, contradicting maximality.
\end{proof}

%
\subsection{Some automorphisms of free group factors and their mixing properties}
Let $\alpha$ be a trace-preserving automorphism of $(L\mb Z,\tau)$, where $\tau(x) = \ip{x\delta_e,\delta_e}$ is the canonical trace. It induces automorphisms of $L\mb F_3$ obtaining by decomposing $L\mb F_3$ as a free product with respect to the subalgebras $(L\mb Z)_a$, $(L\mb Z)_b$, $(L\mb Z)_c$ generated by $a$, $b$ and $c$ using the free product decomposition $L\mb F_3 = (L\mb Z)_a \ast (L\mb Z)_b \ast (L\mb Z)_c$.
We denote
\[
 \Aut_a(L\mb F_3) = \{\alpha\ast \id \in \Aut(L\mb F_3 \cong (L\mb Z)_a\ast L\mb F_2)\,|\,\alpha\in \Aut(L\mb Z,\tau)\},
\]
\[
 \Aut_b(L\mb F_3) = \{\id\ast\beta\ast \id \in \Aut(L\mb F_3 \cong L\mb Z\ast (L\mb Z)_b\ast L\mb Z)\,|\,\beta\in \Aut(L\mb Z,\tau)\}.
\]
\[
 \Aut_c(L\mb F_3) = \{\id\ast\gamma \in \Aut(L\mb F_3 \cong L\mb F_2\ast (L\mb Z)_c)\,|\,\gamma\in \Aut(L\mb Z,\tau)\}.
\]

We get actions of the groups $G_c = \Aut_c(L\mb F_3)$, $G_{bc}= \Aut_b(L\mb F_3)\times \Aut_c(L\mb F_3)$, and $G = \Aut_a(L\mb F_3)\times \Aut_b(L\mb F_3)\times \Aut_c(L\mb F_3)$ on $L\mb F_3$.

The following definition of the relative mixing property appeared in \cite[Definition 2.9]{popa-cocycle-and-oe-superrigidity-malleable}. We will use a strengthening of it.

\begin{Def}
 Let $N\subset M$ be a trace-preserving inclusion of finite von Neumann algebras and let $\sigma\colon G\to \Aut(M)$ be a trace-preserving action of a group $G$ on $M$ such that $\sigma_g(N) = N$ for all $g\in G$. The action $\sigma$ is called weakly mixing relative to $N$ if for every finite set $F\subset M\ominus N$ and for every $\eps > 0$ there is a $g\in G$ such that
 \[
  \norm{E_N(y^*\sigma_g(x))}_2 < \eps,\quad x,y\in F.
 \]
 An action $\sigma$ which is weakly mixing relative to $N=\mb C$ is called weakly mixing.
\end{Def}

When $N$ is pointwise fixed by the action, we get the following result resembling the classical equivalent characterisations of weakly mixing actions (cf. \cite[Proposition D.2]{vaes-rigidity-bourbaki}).

\begin{Prop}\label{prop:weak-mixing-conditions}
 Let $N\subset M$ be a trace-preserving inclusion of finite von Neumann algebras $\sigma\colon G\to \Aut(M)$ be a trace-preserving action such that $\sigma_g(n) = n$ for all $g\in G$ and $n\in N$. Then each of the following conditions implies the next one:
 \begin{enumerate}
  \item $\sigma$ is weakly mixing relative to $N$;
  \item for every $x_1,\dots, x_n \in M\ominus N$ there exists a sequence $g_j\in G$ such that for every $y\in M$ $\norm{E_N(y^*\sigma_{g_j}(x_i))}_2\to 0$, $j\to\infty$, $i=1,\dots,n$.
  \item every finite-dimensional invariant subspace of $M$ is contained in $N$;
  \item for every action $\rho$ of $G$ on a finite von Neumann algebra $P$
  \[
   (M\htens P)^{\sigma\otimes \rho} = N\htens P^\rho.   
  \]
 \end{enumerate}
\end{Prop}
\begin{proof}
The proof of \cite[Prop. D.2]{vaes-rigidity-bourbaki} goes through with small modifications. i)~$\Rightarrow$~ii) is immediate. To show ii)~$\Rightarrow$~iii), take a finite-dimensional invariant subspace $V\subset M$. Consider the space $V':=(1-E_N)(V)\subset M\ominus N$; it is also finite-dimensional and $G$-invariant. But in view of ii) we have for every $x\in V'$ and $y\in M$ that
\[
 \tau(E_N(y^*\sigma_{g_j}(x))) = \ip{\sigma_{g_j}(x),y} \to 0,\quad j\to \infty.
\]
It means that $\sigma_{g_j}(x)\to 0$ weakly as $j\to \infty$. As $V'$ is finite-dimensional, this implies convergence to zero in norm, and as the $G$-action is norm-preserving, this means that $x = 0$. Thus, $V'=\{0\}$, thus $V$ lies in $N$.

To show iii)~$\Rightarrow$~iv), take $T\in (M\htens P)^{\sigma\otimes \rho}\subset  L^2(M)\otimes L^2(P)$ and view it as a Hilbert--Schmidt operator $T\colon \overline{L^2(P)}\to L^2(M)$. Then the image of $T$ is contained in $M$, and the operator $TT^*$ is trace-class and commutes with the $G$-action. Taking its spectral projection, we obtain a finite-dimensional $G$-invariant subspace of $M$, which is necessarily contained in $N$ by iii).  Thus, the image of $T$ lies in $N$, and therefore $T\in N\htens P^\rho$.
\end{proof}

\begin{Lemma}\label{lemma:weak-mixing-free-product}
Let $Q$ and $N$ be finite von Neumann algebras and $\sigma\colon G\to \Aut(Q)$ a trace-preserving weakly mixing action of a group $G$ on $Q$. Then the action $\sigma\ast\id$ of $G$ on the free product von Neumann algebra $M = Q\ast N$ is weakly mixing relative to $N$.
\end{Lemma}
\begin{proof}

By a standard density argument it is enough to check the relative weak mixing condition on the algebraic free product $\mc M = Q\ast_{\mathrm{alg}} N \subset M$, which is the algebra generated by $Q$ and $N$ inside $M$. It is spanned by $N$ and alternating products of elements from $N\ominus\mb C 1$ und $Q\ominus \mb C 1$; without loss of generality we may and will assume that the operator norms of all factors are bounded by 1. Moreover we have
\[
 E_N(n) = n,\quad n\in N
\]
and
\[
 E_N(q_1n_1\cdots q_{k-1}n_{k-1}q_k) = 0
\]
for $q_i \in Q\ominus\mb C 1,\,n_i\in N\ominus\mb C 1$. Thus, for $q_i, q_i'\in Q\ominus\mb C 1,\,n_i,n_i'\in N\ominus\mb C 1$ we get
\[
 E_N(q_1n_1\cdots q_k n_k n'_\ell q'_\ell\cdots n'_1 q'_1) = 0,\quad k\neq \ell,
\]
\[
 E_N(q_1n_1\cdots q_k n_k n'_k q'_k\cdots n'_1 q'_1) = E_N(n_1q_1\cdots n_k q_k q'_k n'_k\cdots q'_1 n'_1) = \prod_{i=1}^k \tau(n_in'_i) \tau(q_i q'_i).
\]

Thus, $\mc M\ominus \mc N$ is spanned by the alternating products of elements from $N\ominus\mb C 1$ and $Q\ominus \mb C 1$, and it's enough to check the weak mixing property for such alternating products. Given a finite set $F$ of them, let $F'\in Q\ominus \mb C 1$ be all factors from $Q\ominus \mb C1$ occuring in the products. By the weak mixing property we find a $g\in G$ such that
\[
 |\tau(q\sigma_g(q'))| < \eps < 1,\quad q,q'\in F'.
\]
The above formulae for $E_N$ then imply for two elements $x,y\in F$ that $E_N(y^* \sigma_g(x)) = 0$ unless
\[
 y^* = n q_1n_1\cdots q_k n_k
\]
\[
 x = n'_k q'_k\cdots n'_1 q'_1 n',
\]
where $q_i, q_i'\in Q\ominus\mb C 1,\,n_i,n_i'\in N\ominus\mb C 1$, $\norm{n_i}\leqslant 1$, $\norm{n'_i}\leqslant 1$, $\norm{q_i}\leqslant 1$, $\norm{q'_i}\leqslant 1$, $n,n'\in N$ with $\norm{n}\leqslant 1$, $\norm{n'}\leqslant 1$.

In this case we get
\[
 \norm{E_N(nq_1n_1\cdots q_k n_k (\sigma_g\ast \id)(n'_k q'_k\cdots n'_1 q'_1n'))}_2 \leqslant \prod_{i=1}^k |\tau(n_in'_i) \tau(q_i \sigma_g(q'_i))| < \eps
\]
which proves the statement.
\end{proof}

Using the existence of weakly mixing actions on an arbitrary finite von Neumann algebra (e.g. Bernoulli actions \cite[Sect. 2.4]{popa-nc-bernoulli}) and Proposition \ref{prop:weak-mixing-conditions}, we obtain
\begin{Cor}\label{lemma:perm-invariant}
Consider the actions of $G_c = \Aut_c(L\mb F_3)$, $G_{bc}= \Aut_b(L\mb F_3)\times \Aut_c(L\mb F_3)$, and $G = \Aut_a(L\mb F_3)\times \Aut_b(L\mb F_3)\times \Aut_c(L\mb F_3)$ on $L\mb F_3$ described above. Then
\begin{enumerate}
 \item every $G_{bc}$-invariant element $L\mb F_3\htens L\mb F_3^\op$ is contained in $(L\mb Z)_a\htens (L\mb Z)_a^\op$;
 \item every $G_{c}$-invariant element $L\mb F_3\htens L\mb F_3^\op$ is contained in $(L\mb F_2)_{ab}\htens (L\mb F_2)_{ab}^\op$.
\end{enumerate}
\end{Cor}

\section{Non-continuous derivations}
Let $\mb F_3 = \langle a,b,c \rangle$ be a free group on three generators. We will naturally view $\mb F_2 = \langle b,c \rangle$ as a subgroup of $\mb F_3$. It's well-known that $\mb F_2=\langle b,c\rangle$ contains a copy of $\mb F_\infty = \langle g_1,\,g_2,\dots \rangle$, and we well fix such a copy.

We recall that the derivation
\[
 \partial_a\colon \mb C[\mb F_3]\to \ms U(L\mb F_3\htens L\mb F_3^\op)
\]
is uniquely defined by the conditions
\[
 \partial_a(a) = a\otimes 1^\op,\quad \partial_a(b)=\partial_a(c) = 0.
\]
\begin{Prop}\label{prop:no-continuous-extension}
Let $p \in L\mb F_3\htens L\mb F_3^\op$ be a nonzero projection. There is no norm-measure continuous extension of $\partial_a\cdot p$ to $L\mb F_3$.
\end{Prop}
\begin{proof}
 We will first construct a particular sequence $y_n \in \mb C[\mb F_3]$ such that $\norm{y_n}_\infty \to 0,\,n\to\infty$, but $\partial_a(y_n)\not\xrightarrow{m} 0$. We set
 \[
  x_n := g_1ag_2ag_3a\cdots ag_na\in \mb C[\mb F_3]
 \]
and consider
\[
 \partial_a(x_n) = \sum_{k=1}^{n} g_1a\cdots g_ka \otimes g_{k+1}a\cdots g_na
\]
The elements
\[
 g_1a,\, g_1ag_2a,\dots,\, g_1ag_2a\cdots g_na \in \mb F_3
\]
are free. Indeed, the family $\{g_ia\}_{i=1}^n$ is free, because $\{g_i\}_{i=1}^n$ is a free family which is itself free from $a$.

Therefore the elements
\[
 (g_1a\cdots g_k a, g_{k+1}a\cdots g_na) \in \mb F_3\times \mb F_3^\op,\quad k=\overline{1,n}
\]
are free and hence form a freely independent family of unitaries $\{u_i\}_{i=1}^n\subset L\mb F_3\htens L\mb F_3^\op$. Consider the elements
\[
 h_n:=\Re \partial_a(x_n) = \frac{1}{2}(\partial_a(x_n) + \partial_a(x_n)^*) = \frac{1}{2}\sum_{i=1}^n (u_i+u_i^*) \in L\mb F_3\htens L\mb F_3^\op.
\]
%
The spectral density of this operator can be computed using free probability theory. Indeed, $h_n$ is an instance of a scaled random walk operator on a free group, and therefore we can use the formula for the Cauchy transform for its spectral density from \cite[Example 3.4.5]{voiculescu-dykema-nica}:
\[
 G_{h_n}(\zeta) = \frac{n\zeta\sqrt{1-(2n-1)\zeta^{-2}}-(n-1)}{\zeta^2 - n^2}
\]
Using the Stieltjes inversion formula, we get the spectral density of $h_n$:
\[
 d\mu_{h_n}= \begin{cases}
              \dfrac{n\sqrt{(2n-1)-x^2}}{n^2-x^2}dx,& x\in[-\sqrt{2n-1},\sqrt{2n-1}],\\
              0 & \text{otherwise}.
             \end{cases}
\]


We see that
\[
 \tau\left(\chi_{[-\eps,\eps]}\left(\frac{h_n}{\sqrt[4]{n}}\right)\right) = \int_{-\eps\sqrt[4]{n}}^{\eps\sqrt[4]{n}} \frac{n\sqrt{(2n-1)-x^2}\,dx}{n^2-x^2}\leqslant \frac{2\eps\cdot n^{5/4}\sqrt{2n-1}}{n^2-\eps^2\sqrt{n}}\to 0,\quad n\to\infty.
\]
In view of Lemma \ref{lemma:no-convergence-to-zero}, for every nonzero projection $p\in L\mb F_3\htens L\mb F_3^\op$ we get
\[
 \frac{h_n}{\sqrt[4]{n}} \cdot p \not\xrightarrow{m} 0,\quad n\to\infty.
\]

In particular,
\[
 \partial_a\left(\frac{x_n}{\sqrt[4]{n}}\right) \not\xrightarrow{m} 0,\quad n\to\infty,
\]
although $x_n/\sqrt[4]{n}$ converges to $0$ in norm, because $\norm{x_n}_\infty = 1$. Thus, $y_n:= x_n/\sqrt[4]{n}$ satisfies the required properties, and the derivation $\partial_a$ is not continuous.

Now we have to show that the continuity projection $q$ of $\delta$ is equal to $0$. The derivation $\partial_a$ is invariant under the actions of the groups $\Aut_b(L\mb F_3)$ and $\Aut_c(L\mb F_3)$. By Lemma \ref{lemma:invariant-projection} we get that $q$ is invariant under $\Aut_b(L\mb F_3)$ and $\Aut_c(L\mb F_3)$, hence by Corollary \ref{lemma:perm-invariant}, i) it belongs to $(L\mb Z)_a\htens (L\mb Z)_a^\op$.

Now, for $\alpha\in \Aut_a(L\mb F_3)$ we infer that
\[
 \partial^\alpha_a(x) = \alpha^{-1}(\partial_a(\alpha(x)))\alpha^{-1}(q)=\partial_a\cdot (a^{-1}\otimes 1^\op)\alpha^{-1}(\partial_a(\alpha(a)))\alpha^{-1}(q)
\]
is continuous. The support projection of $(a^{-1}\otimes 1^\op)\alpha^{-1}(\partial_a(\alpha(a)))\alpha^{-1}(q)$ is equal to $\alpha^{-1}(q)$, and by maximality of $q$ it follows that $\alpha^{-1}(q) = q$. Therefore, $q=0$ or $q=1$, but the latter case is impossible because $\partial_a$ is not continuous.

\end{proof}

For obvious symmetry reasons, the statement of Proposition \ref{prop:no-continuous-extension} is also true for the derivations
\[
 \partial_b,\,\partial_c\colon \mb C[\mb F_3]\to \ms U(L\mb F_3\htens L\mb F_3^\op)
\]
uniquely determined by the conditions
\[
 \partial_b(b) = b\otimes 1^\op,\quad \partial_b(a)=\partial_b(c) = 0
\]
resp.
\[
 \partial_c(c) = c\otimes 1^\op,\quad \partial_c(a)=\partial_c(b) = 0.
\]


Is is well-known that for every nonzero projection $p\in L\mb F_3\htens L\mb F_3^\op$ the derivations $\partial_a\cdot p$, $\partial_b\cdot p$, $\partial_c \cdot p$ are not inner. The derivations $\partial_a$, $\partial_b$, $\partial_c$ freely generate the module $\Der(\mb C\mb F_3,\ms U(L\mb F_3\htens L\mb F_3^\op))$.

\begin{Thm}
 The first continuous $L^2$-cohomology of $L\mb F_3$ vanishes:
 \[
 H^1_c(L\mb F_3,\ms U(L\mb F_3\htens L\mb F_3^\op))\cong 0.
 \]
\end{Thm}
\begin{proof}
The restriction of $\delta$ to $\mb C\mb F_3$ is a derivation and therefore can be uniquely written as a combination of $\partial_a$, $\partial_b$, $\partial_c$:
 \[
  \delta(x) = \partial_a(x)\cdot \xi''_a + \partial_b(x)\cdot \xi''_b + \partial_c(x)\cdot \xi''_c,\quad x\in\mb C\mb F_3
 \]
for some $\xi''_a,\,\xi''_b,\,\xi''_c\in \ms U(L\mb F_3\htens L\mb F_3^\op)$. As inner derivations are continuous, after subtracting an inner derivation we may assume that $\xi_c = 0$ and
\[
 \delta(x) = \partial_a(x)\cdot \xi'_a + \partial_b(x)\cdot \xi'_b,\quad x\in\mb C\mb F_3.
\]

The right $\ms U$-module generated by $\partial_a\cdot \xi_a$ contains $\partial_a\cdot p_a$, where $p_a$ is the target projection of $\xi_a$. As the right $\ms U$-action preserves continuity, we may assume that $\xi_a=p_a$ and that $\delta$ has the form
\[
 \delta(x) = \partial_a(x)\cdot p_a + \partial_b(x)\cdot \xi_b,\quad x\in\mb C\mb F_3.
\]
Multiplying from the right with $(1-p_a)$ and using Proposition \ref{prop:no-continuous-extension}, we deduce $\xi_b(1-p_a) = 0$. Thus, $\rk \xi_b \leqslant \rk p_a$; reasoning symmetrically, we infer $\rk \xi_b = \rk \xi_a = \rk p_a$. We also observe that for every $\delta \in \Der_c(\mb C\mb F_3,\ms U(L\mb F_3\htens L\mb F_3^\op))$, the elements $p_a$ and $\xi_b$ are uniquely determined.


Now, let
\[
 P_{\rm{cont}} = \{p\in \Proj(L\mb F_3\htens L\mb F_3^\op)\,|\, \exists \delta \in \Der_c(\mb C\mb F_3,\ms U(L\mb F_3\htens L\mb F_3^\op))\colon p=p_a\}.
\]
Analogously to Lemma \ref{lemma:complete-lattice}, we are going to prove that $P_{\rm{cont}}$ is a complete lattice. Indeed, if $p_1,p_2\in P_{\rm{cont}}$, then $p_1\wedge p_2 \in P_{\rm{cont}}$ because $p_1\cdot (p_1\wedge p_2) = p_1\wedge p_2$. Now, if $\delta_1$ and $\delta_2$ are derivations corresponding to $p_1$ resp. $p_2$, then $p_1\vee p_2$, being the support projection of $p_1+p_2$, corresponds to the derivation $\delta_1+\delta_2$. Completeness of $P_{\rm{cont}}$ is proven as follows: if let $(p_i)_{i\in I}$ be an orthogonal family in $P_{\rm{cont}}$ with corresponding derivations $\delta_i$ with elements $\xi'_{b,i}$. By the observation above, $\rk(\xi'_{b,i}) = \rk(p_i)$. But then
\[
 \sum_{i\in I} \rk(\xi'_{b,i}) = \sum_{i\in I} \rk(p_i) \leqslant 1
\]
and therefore as in Lemma \ref{lemma:complete-lattice} the series
\[
 \sum_{i\in I} \delta_i(x) = \sum_{i\in I} (\partial_a(x)\cdot p_i + \partial_b(x)\cdot \xi_{b,i})
\]
converges uniformly in measure to a derivation $\delta$ having the supremum of $p_i$ as the corresponding projection $p_a$.

Thus, $P_{\rm{cont}}$ is a complete lattice. Let $p$ be its maximal element. In view of the equalities
\[
 \delta^\beta(x) =  \partial_a(x)\cdot \beta^{-1}(p_a) + \beta^{-1}(\partial_b(\beta(x)))\beta^{-1}(\xi_b)
\]
and
\[
 \delta^\gamma(x) = \partial_a(x)\cdot\gamma^{-1}(p_a) + \partial_b(x)\cdot \gamma^{-1}(\xi_b)
\]
for $\beta \in \Aut_b(L\mb F_3)$ and $\gamma\in\Aut_c(L\mb F_3)$ we get by maximality that $p$ is invariant under $\Aut_{b}(L\mb F_3)$ and $\Aut_c(L\mb F_3)$. Therefore by Corollary \ref{lemma:perm-invariant},~i) we deduce that $p\in (L\mb Z)_a\htens (L\mb Z)_a^\op$.


Now, for $\alpha\in \Aut_a(L\mb F_3)$ we get that
\begin{multline*}
 \delta^\alpha(x) = \alpha^{-1}(\partial_b(\alpha(x)))\alpha^{-1}(p) + \partial_b(x)\cdot \alpha^{-1}(\xi_b)\\= \partial_a\cdot (a^{-1}\otimes 1^\op)\alpha^{-1}(\partial_a(\alpha(x)))\alpha^{-1}(p) + \partial_b(x)\cdot \alpha^{-1}(\xi_b) 
\end{multline*}
is continuous. The support projection of $(a^{-1}\otimes 1^\op)\alpha^{-1}(\partial_a(\alpha(a)))\alpha^{-1}(p)$ is equal to $\alpha^{-1}(p)$, and by maximality of $p$ it follows that $\alpha^{-1}(p) = p$. Therefore, $p=0$ or $p=1$.

If $p=1$, we are given a continuous derivation of the form
\[
  \delta_0(x) = \partial_a(x) + \partial_b(x) \cdot \zeta_b,\quad x\in\mb C\mb F_3.
\]
For $\gamma\in \Aut_c(L\mb F_3)$ we obtain
\[
  \delta_0^\gamma(x) = \partial_a(x) + \partial_b(x)\cdot \gamma^{-1}(\zeta_b),\quad x\in\mb C\mb F_3.
\]
As $\delta_0^\gamma - \delta_0$ is continuous, Proposition \ref{prop:no-continuous-extension} implies that
\[
 \gamma^{-1}(\zeta_b) = \zeta_b,\quad \gamma\in \Aut_c(L\mb F_3).
\]
Thus using Corollary \ref{lemma:perm-invariant},~ii), we get that $\zeta_b \in \ms U(L\mb F_2\htens L\mb F_2^\op)$, where $L\mb F_2\subset L\mb F_3$ is generated by $a$ and $b$. Consider the restriction of the derivation $\delta$ to $L\mb F_2$. Subtracting an inner derivation and multiplying with a suitable element of $\ms U(L\mb F_2\htens L\mb F_2^\op)$ from the right, we may assume that the derivation
\[
  \delta_1(x) = \partial_a(x)r_a,\quad x\in\mb C\mb F_2,
\]
is continuous for some nonzero projection $r_a$. Thus, the continuity projection of $\partial_a\colon \mb C\mb F_2 \to \ms U(L\mb F_2\htens L\mb F_2^\op)$ is nonzero. Arguing as in Proposition \ref{prop:no-continuous-extension}, we deduce that the continuity projection of $\partial_a\colon \mb C\mb F_2 \to \ms U(L\mb F_2\htens L\mb F_2^\op)$ is equal to 1. For symmetry reasons, $\partial_b\colon \mb C\mb F_2 \to \ms U(L\mb F_2\htens L\mb F_2^\op)$ is continuous as well. Therefore the module of continuous derivations $\Der_c(L\mb F_2,\ms U(L\mb F_2\htens L\mb F_2^\op))$ is two-dimensional, and hence the first continuous $L^2$-Betti number of $L\mb F_2$ is equal to $1$:
\[
 \eta_1^{(2)}(L\mb F_2) = 1.
\]
But then by the compression formula (Theorem \ref{thm:compression-formula}) we obtain
\[
 \eta_1^{(2)}(L\mb F_3) = 2, 
\]
which contradicts the result of Proposition \ref{prop:no-continuous-extension} that the submodule generated by $\partial_a$ consists of discontinuous derivations.
\end{proof}

As a corollary we get the following result.
\begin{Thm}
 The first continuous $L^2$-cohomology of an interpolated free group factor $\mb F_r$, $1<r<\infty$, vanishes:
 \[
 H^1_c(L\mb F_r,\ms U(L\mb F_r\htens L\mb F_r^\op))\cong 0.
 \]
\end{Thm}
\begin{proof}
 This follows immediately from the case $r=3$ by the compression formula (Theorem \ref{thm:compression-formula}) and Proposition \ref{prop:vanishing-cohomology-vs-vanishing-betti}.
\end{proof}

In particular, the answer to the question of Andreas Thom in \cite{thom2008l2} is negative: Voiculescu's free difference quotients don't have continuous extensions to $L\mb F_n$. Our result also allows to ask the following question.

\subsection*{Question} Is $\eta_1^{(2)}(M) = 0$ for all II$_1$-factors $M$?

\end{document}